\documentclass[centertags,leqno,11pt]{article}

\usepackage{amsmath}
\usepackage{amssymb}
\usepackage{latexsym}
\usepackage{amsxtra}
\usepackage{amscd}
\usepackage{theorem}

\usepackage{graphics}
\usepackage{epic}

\theoremstyle{change}

{\theorembodyfont{\slshape}

\newtheorem{thm}{Theorem.}[section]

\newtheorem{prop}[thm]{Proposition.}

\newtheorem{defn}[thm]{Definition.}}

{\theorembodyfont{\rmfamily}

\newtheorem{rem}[thm]{Remark.}

}

\renewcommand{\em}{\sl}

\newcommand{\proof}{\noindent {\bf Proof:\ }}
\newcommand{\Endproof}{\hspace*{\fill} $\Box$ \vspace{1ex} \noindent }

\makeatletter
\renewcommand{\subsection}{\@startsection{subsection}{2}%
{\z@}{-3.25ex plus -1ex minus-.2ex}{-1em}{\bf}} \makeatother

\newcommand{\PP}{\mathbb{P}}

\newcommand{\CC}{\mathbb{C}}

\renewcommand{\AA}{\mathbb{A}}
\newcommand{\GG}{\mathbb{G}}

\newcommand{\h}{{\bf h}}
\newcommand{\HHH}{{\cal H}}

\newcommand{\x}{{\bf x}}
\newcommand{\g}{{\bf g}}

\newcommand{\GL}{{\rm GL}}
\newcommand{\SL}{{\rm SL}}

\newcommand{\Hom}{{\rm Hom}}

\newcommand{\MC}{{\rm MC}}

\newcommand{\et}{{\rm\acute{e}t}}

\newcommand{\LS}{{\rm LS}}

\newcommand{\codim}{{\rm codim}}
\newcommand{\Sp}{{\rm Sp}}

\newcommand{\To}{\;\longrightarrow\;}

\newcommand{\Mapsto}{\;\longmapsto\;}

\newcommand{\diag}{{\rm diag}}

\newcommand{\pr}{{{\rm pr}}}

\newcommand{\J}{{\rm J}}

\newcommand{\GGG}{{\cal G}}

\newcommand{\SO}{{\rm SO}}

\renewcommand{\L}{{\cal L}}

\numberwithin{equation}{subsection}
\numberwithin{thm}{subsection}
\theoremstyle{plain}

\title{On exceptional rigid local systems}

\author{
   Michael Dettweiler  and 
 Stefan Reiter
  }

\begin{document}
\maketitle 
\begin{abstract} We prove new instances 
of Simpson's rigidity conjecture which states that 
quasi-unipotent rigid local systems should be  motivic. 
We construct new relative motives over the
fourfold punctured Riemann sphere which give rise
to  $G_2$-rigid local systems which are not 
rigid in the group  $\GL_7.$

\end{abstract}


\section*{Introduction}\label{Introduction}

Let $G$ be a reductive algebraic group over $\CC$ and 
let $X$ be  a smooth quasi-projective complex variety.
Let us call a representation $\rho:\pi_1(X)\to G$ to be 
{\it $G$-rigid}, if the set
theoretic orbit of the representation $\rho$ in the representation
space $\Hom(\pi_1(X),G)$ under the action of $G$ is an open subset.
 
Note that by fixing an embedding of $G$ into a  general linear group
$\GL_n(\CC),$ any 
representation
$\rho :\pi_1(X)\to G$ corresponds uniquely  to  a local system $\L_\rho$ 
on $X,$
see \cite{Deligne70}. On the other hand, any local system on $X$ 
whose monodromy lies
in $G\leq \GL_n$
 gives rise to a representation  $\rho=\rho_\L :\pi_1(X)\to G.$
The local system $\L=\L_\rho$ is called {\it $G$-rigid}, if 
$\rho$ is $G$-rigid. A local system is called 
{\it quasi-unipotent}, if the 
eigenvalues of the local monodromies of $\L$ are roots of unity. 

The following conjecture is motivated
by  Hodge theory and appears in 
work of Simpson \cite{Simpson92} (for projective $X$):

\vspace{.3cm}

\noindent {\bf Rigidity Conjecture:} {\rm Let $G$ be a reductive 
complex algebraic group and 
let $X$ be  a smooth quasi-projective complex variety. Then any 
quasi-unipotent 
$G$-rigid local system on $X$ 
is {\it motivic}, i.e., it is a subfactor of  a higher direct image sheaf
$R^i\pi_*(\CC),$ where $\pi:Y\to X$ is some smooth morphism.}\\

The case of $\GL_n$-rigidity for 
representation of $\pi_1(X),$ where $X$ is the punctured Riemann sphere, 
is probably the most important one.
Many classical ordinary 
differential equations (e.g., the 
generalized hypergeometric differential equations)
give rise to $\GL_n$-rigid local systems on
$\PP^1\setminus \{x_1,\ldots,x_r\}$ -- and can be studied in detail using the concept of 
rigidity, see e.g. \cite{Beukers-Heckman} and \cite{Katz96}.  Moreover, by the
work of N. Katz \cite{Katz96}, there exist a 
uniform description of all irreducible $\GL_n$-rigid local systems 
(loc.~cit., Chap. 6) 
 and the 
rigidity conjecture has been proven for irreducible
$\GL_n$-rigid local systems (loc. cit., Chap. 8). \\

Let $G\leq H$ denote an embedding of algebraic groups (think about
 the embedding
of the exceptional simple group $G_2$ into the general linear group 
$\GL_7$). If $\L$ is a $G$-rigid local system, then in most cases $\L$ is not 
$H$-rigid. 
In a recent paper \cite{DR06}, we use Katz' work 
and certain $G_2$-rigid representations of $\pi_1(\PP^1\setminus
\{0,1,\infty\})$
which are (surprisingly)  $\GL_7$-rigid, in order to construct
motivic $G_2$ local systems. 
Unfortunately, most $G_2$-rigid local systems are not 
$\GL_7$-rigid -- and the rigidity conjecture for them 
cannot be deduced from
Katz' work so easily. \\

It is the aim of this article to show that nevertheless,  some  
$G_2$-rigid local systems are motivic although 
they are not $\GL_7$-rigid. We prove the following result
(where $\J(n_1,\ldots,n_k)$ denotes a unipotent matrix in Jordan 
canonical form which decomposes into blocks of length $n_1,\ldots,n_k$
and $\zeta_3$ denotes a primitive third root of unity):\\

\noindent {\bf Theorem 1:} {\it Let  $x_1,x_2,x_3$ be pairwise distinct
complex numbers. Then the following holds:
\begin{enumerate}
\item There exists a $G_2(\CC)$-rigid
local system $\HHH$  on $\PP^1\setminus \{x_1,x_2,x_3,\infty\}$
such that the Jordan canonical forms of the 
local monodromies at $x_1,x_2,x_3,\infty$ are as follows:
$$ \J(2,2,1,1,1),\,\,\J(2,2,1,1,1),\, \,
\diag(1,\zeta_3,\zeta_3,\zeta_3,\zeta_3^{-1},\zeta_3^{-1},\zeta_3^{-1}),\,\, 
\J(3,3,1)\,.$$
\item The Zariski closure of the image of the monodromy of $\HHH$ 
coincides with the exceptional simple group $G_2(\CC).$
\item 
The local system $\HHH$  is $\SO_7(\CC)$-rigid but not $\GL_7(\CC)$-rigid.
\item 
The rigidity conjecture holds for $\HHH.$  
\end{enumerate}}

The proof of Theorem 1 is given in Section \ref{Sec5}. An outline of the 
proof is as follows: We start with two $\GL_2(\CC)$-rigid local systems
 and take an 
appropriate tensor product  of them. This gives rise to a local system 
of rank four on $\PP^1\setminus \{x_1,x_2,x_3,\infty\}$
which is $\SO_4(\CC)$-rigid
but not $\GL_4(\CC)$-rigid. 
 Then we apply the middle convolution and suitable tensor 
operations twice
end up with  the above $G_2(\CC)$-rigid local system $\HHH$ (it is constructed 
in Section \ref{Sec5}).
The rigidity conjecture for $\HHH$
follows then from the 
motivic interpretation of the middle convolution as given in the work of Katz \cite{Katz96}, Chap. 8. 
It seems remarkable, 
that the weight 
of $\HHH$ (as a variation of Hodge-structures)
is $4$ contrary to the $G_2$-local system  
considered in \cite{DR06} which has weight $6.$\\

This article has grown out from a question of P. Deligne about 
a possible motivic interpretation of the 
other $G_2$-rigid (but non-$\GL_7$-rigid) local systems that exist 
besides the $G_2$-rigid local systems 
 considered in \cite{DR06}
and some related $G_2$-rigid local systems  which where 
communicated to the authors 
by P. Deligne and N. Katz. One may hope that the above 
approach leads to a verification of the rigidity conjecture for 
rigid local systems with values in other exceptional groups.

The authors thank P. Deligne and N. Katz for their interest in the 
local system of \cite{DR06} and for several valuable remarks on possible
generalizations of the results in loc.~cit..

\section{Preliminaries} \label{Sec4}

\subsection{Rigid local systems and rigid tuples.}\label{Sec41}
Let $X$ be $\PP^1\setminus D\,,$ where
$D=\{x_1,\ldots, x_r\}\subseteq \PP^1.$ Fix 
a homotopy base 
$\gamma_1,\ldots,\gamma_r$ of $\pi_1(X)$ satisfying the product relation 
$\gamma_1\cdots \gamma_r=1.$ 
The category of 
local systems of $\CC$-vector spaces will be denoted by $\LS(X).$  
It is well known that any local system corresponds uniquely
to its {\it monodromy representation} 
$$ \rho_\L: \pi_1(X,x_0)\To \GL(\L_{x_0})\,,
$$ see e.g. \cite{Deligne70}. 
Thus any local system $\L\in \LS(X)$ of rank $n$
corresponds to its {\it monodromy  tuple} $\g_{\L}=(g_1,\ldots,g_r),$ where 
$g_i \in \GL(\L_{x_0})\simeq \GL_n(\CC)$ is the image of $\gamma_i$ under 
the monodromy representation  $\rho_\L.$ By construction, the 
product relation $g_1\cdots g_r=1$ holds for any monodromy tuple 
$\g_\L,\, \L\in \LS(X).$ (The concept of a 
monodromy tuple enables in 
many cases the explicit computation of a local system
$\L\in \LS(X).$)\\

Let $G\leq \GL_n(\CC)$ be a reductive complex algebraic group.
By definition, a  local system $\L\in \LS(X)$ 
with monodromy tuple $\g=(g_1,\ldots,g_r)\in G^r$ (resp. $\rho_\L$)
is {\it $G$-rigid}, if  there are (up to simultaneous $G$-conjugation) 
only finitely many tuples $\h=(h_1,\ldots,h_r)\in G^r$ with 
$ h_1\cdots h_r=1$ and such that $h_i$ is $G$-conjugate to $g_i.$ In this 
case  the tuple $\g$ is called {\it $G$-rigid}.\\

The following theorem is often useful to detect
$G$-rigid tuples: 

\begin{thm} \label{propSV} Let $G\leq \GL_n(\CC)$ be a reductive algebraic 
subgroup. Let $\L$ be an irreducible local system of rank $n$ 
whose monodromy tuple  
$\g=(g_1,\ldots,g_r)$ is contained in 
$ G^r.$ 
Then the local system $\L$ (resp. its monodromy  tuple $\g$)
 is 
$G$-rigid, if and only if the following {\rm dimension formula} holds:
$$ \sum_{i=1}^r \codim(C_G(g_i))=2(\dim(G)-\dim(Z(G)),$$
where $C_G(g_i)$ denotes the 
centralizer of $g_i$ in $G,$
the codimension is taken relative to $G,$ and $Z(G)$ denotes the centre 
of $G.$
\end{thm}

\proof The necessity of the dimension formula for $G$-rigidity
is proven in \cite{SV}.
The other direction follows from the same arguments
as in \cite{Weil64}, Section 7.
\Endproof

\subsection{Operations on local systems}\label{Sec42}
The following operations on local systems on $X$ will play a role
below:\\

On the one hand side, there is the usual tensor product
$\L_1\otimes \L_2$  of 
local systems $\L_1,\L_2\in \LS(X)$ having ranks 
$n_1,n_2,$ respectively. 
Let $\g_{\L_1}=(g_1,\ldots,g_r)$ and $\g_{\L_2}=(h_1,\ldots h_r)$ be 
the monodromy tuples of $\L_1$ and $\L_2$ (resp.). Then 
$$ \g_{\L_1\otimes \L_2}=\g_{\L_1}\otimes \g_{\L_2}:=(g_1\otimes h_1,\ldots, g_r\otimes h_r)\,\in \, 
\GL_{n_1n_2}^r,$$ 
where $g_i\otimes h_i$ denotes the usual 
Kronecker product of matrices. Thus, 
the tensor product of local systems is very well understood from the 
computational side. Moreover, if $\L_1$ and $\L_2$ are irreducible 
and motivic, then 
it follows from the K\"unneth formula that $\L_1\otimes \L_2$ is also 
motivic. \\

The next operation on local systems on $$X=\AA^1\setminus \x,\, 
\x=\{x_1,\ldots,x_r\}=\PP^1\setminus (\x\cup \infty),$$ was introduced by 
N. Katz in \cite{Katz96} and is much less obvious:  
 Let $p(x)=\prod_{i=1}^r(x-x_i)$ and let 
$L$ be the divisor on $\AA^2=\AA^1_x\times \AA^1_y$
 associated to $p(x)\cdot p(y)\cdot (y-x)=0.$ Let 
$U=\AA^2\setminus  L$ and let 
$$ j: U \To \PP^1_{\AA^1_y\setminus \x} ,\, (x,y) \Mapsto ([x,1],y)\,$$
(by abuse of notation, the divisor on $\AA^1_y$ which is given 
by $p(y)=0$ is also denoted by $\x$). 
Let further $\pr_i:U \to \AA^1\setminus \x,\, i=1,2,$ 
denote the $i$-th projection, let 
$$ {\rm d}:U\To \GG_m=\AA^1\setminus \{0\},\, (x,y)\Mapsto y-x$$ be the {\em difference map} and let $\bar{\pr}_2:\PP^1_{\AA^1\setminus \x} \to {\AA^1\setminus \x}$ be the natural projection. Following Katz \cite{Katz96}, Chap. 8,
we define a middle convolution operation on 
$\LS(\AA^1\setminus \x)$ as follows:

\begin{defn} {\em For a character 
$\chi:\pi_1^\et(\GG_m)\to \CC^\times,$ let $\L_\chi\in \LS(\GG_m)$ denote
the associated local system. For $\L\in \LS(\AA^1\setminus \x),$ the local system
$$ \MC_\chi(\L):=R^1(\bar{\pr}_2)_*\left(
j_*(\pr_1^*(\L)\otimes {\rm d}^*(\L_\chi))\right)\,\in \LS_R(\AA^1\setminus\x)$$
is called  the 
{\em middle convolution of $\L$ with $\L_\chi.$}}
\end{defn}

\begin{rem}\label{rem1}{\rm The above definition amounts to a special case
of the middle convolution of perverse sheaves, introduced 
in \cite{Katz96}.
It follows from \cite{DeHabil}, L. 3.5.5, that if $\L$ is irreducible 
and has at least two nontrivial local monodromies away from 
$\infty,$ then the monodromy tuple
of $ \MC_\chi(\L)$ is given by the tuple $\MC_\lambda(g_\L),$ where
 $\MC_\lambda$ is the tuple-transformation introduced in \cite{DR00}
(the tuple $\MC_\lambda(\g_\L)$ is constructed from $\g_\L$ by
an explicit receipt using only linear algebra). 
Thus, 
the middle convolution  of local systems with Kummer sheaves
is also well understood from the 
computational side. }\end{rem}

The following Proposition can be proved using the 
the same arguments 
as \cite{Katz96}, Chap. 8:

\begin{prop}\label{propkat} Let 
 $\pi:Y\to \AA^1\setminus \x $ be a smooth proper 
morphism and let  $\L$ be an 
irreducible subsheaf of $R^i\pi_*(\CC)\in \LS(\AA^1\setminus \x)$ 
which has  at least two non-trivial 
local monodromies away from $\infty.$ Assume further that $\L$ has 
weight $k$ (as variation of Hodge structures).  For a natural number $n,$
let $$\alpha_n:\GG_m\to \GG_m,\,\, x\mapsto x^n$$
and let $$\Pi:=\pr_2\circ (\pi\boxtimes \alpha_n),$$ where 
 $\pi\boxtimes \alpha_n$ denotes the fibre product 
of the pullbacks of $\pi$ and $\alpha_n$ to $U.$ Let 
$\chi:\pi_1(\GG_m)\to \CC$ be the character of 
$\pi_1(\GG_m)$ which 
sends a counterclockwise 
generator of $\pi_1(\GG_m)$ to $e^{2\pi i/n}.$
Then 
 $ \MC_\chi(\L)\in \LS(\AA^1\setminus \x)$ is  a subfactor of 
$R^1\Pi_*(\CC)\in \LS(\AA^1\setminus \x)$ 
and  the weight of $ \MC_\chi(\L)$ as a variation 
of Hodge structures is equal to $k+1.$
\end{prop}

\subsection{Construction of local systems with finite monodromy}\label{Sec43}
The following construction of local systems
with finite monodromy will be used in Section 
\ref{Sec5} below: Let $\rho_f:\pi_1^\et(X)\to G$ be 
the surjective homomorphism associated to a finite 
Galois cover of $X$ with Galois group $G$ 
and let $\alpha:G\to \GL_n(\CC)$ be a representation.
Then the composition 
$ \alpha\circ \rho_f: \pi_1^\et(X)\to \GL_n(R)$ corresponds to a 
local system 
 with finite monodromy on $X$ which is denoted by $\L_{f,\alpha}.$ Note
that $\L_{f,\alpha}$ is motivic by construction.

\section{The proof of Theorem~1}\label{Sec5}

Let $x_1,x_2,x_3\in \CC$ be pairwise distinct,  let 
$X=\PP^1\setminus \{x_1,x_2,x_3,\infty\}$ and let $\zeta_3$ denote a 
fixed primitive third  root of unity.
 Fix a homotopy base $\gamma_i,\,i=1,2,3,4$
of $\pi_1(X,x_0)$ such that $\gamma_i,\,i=1,2,3,$ is a simple loop around 
the missing point $x_i,$ such that $\gamma_4$ is a
 simple loop around $\infty,$
and such that $\gamma_1\cdots \gamma_4=1.$

\subsection{Construction of the underlying local system}
Using the notation of Section \ref{Sec43}, let $\L_i=\L_{f_i,\alpha}\in 
\LS(X),\, 
i=1,2,3,$ be as follows:

\begin{itemize}
\item
The Galois cover $f_1:Y_1\to X$ is the cover of $X$ 
with Galois group $$Z_3=\langle \sigma\mid \sigma^3=1 \rangle$$ which 
is ramified
at $x_1,x_3$ and $\infty,$ and $\alpha$ sends $\sigma$ to 
$\zeta_3\in \CC^\times=\GL_1(\CC).$ Thus the monodromy tuple of $\L_1\in \LS(X)$ 
is equal to $\g_1=(\zeta_3,1,\zeta_3,\zeta_3).$ 
\item The Galois cover $f_2:Y_2\to X$ is the cover of $X$ 
with Galois group $Z_3$ 
which is ramified
at $x_2,x_3$ and $\infty.$ Thus the monodromy tuple of $\L_2\in \LS(X)$ 
is equal to $\g_2=(1,\zeta_3,\zeta_3,\zeta_3).$ 
\item The Galois cover $f_3:Y_3\to X$ is the cover of $X$ 
with Galois group $Z_3$
which is ramified
at $x_3$ and $\infty.$ Thus  
the  monodromy tuple of $\L_3$ is $\g_3=(1,1,\zeta_3,\zeta_3^{-1}).$
\end{itemize}

Let $\chi:\pi_1(\GG_m)\to \CC^\times$ be the character which sends 
a counterclockwise generator of $\pi_1(\GG_m)$ to $\zeta_3$ and 
let $\chi^{-1}$ be the dual character.
Consider the following sequence of tensor operations and 
middle convolutions of local systems:
\begin{equation}\label{eq1}
\HHH:=\L_3^{-1}\otimes( \MC_{\chi^{-1}}(\L_3\otimes( \MC_{\chi}(\MC_{\chi^{-1}}(\L_1)\otimes \MC_{\chi^{-1}}(\L_2)))))\,.\end{equation}
By Remark \ref{rem1}, The monodromy tuple of $\GGG$ is given by 
\begin{equation}\label{eq2}
\h=(h_1,h_2,h_3,h_4)=\g_3^{-1} \otimes( \MC_{\zeta_3^{-1}}(\g_3\otimes( \MC_{\zeta_3}(
\MC_{\zeta_3^{-1}}(\g_1)\otimes \MC_{\zeta_3^{-1}}(\g_2)))))\,,\end{equation}
where $\g_3^{-1}=(1,1,\zeta_3^{-1},\zeta_3)$ and where we have used 
the convention of Section \ref{Sec42} for the tensor product of tuples. 
Using the explicit receipt for the computation of $\MC_\lambda,$ 
one can evaluate Equation \eqref{eq2} and finds matrices 
for the monodromy tuple 
 $\h$ (these are given in the Appendix).  

\subsection{Proof of Theorem~1:}
The following result gives a proof of Theorem 1, (i)
and (ii) (see Introduction):

\begin{prop} \label{Properg1}  Let 
$$\rho =\pi_1(\PP^1\setminus \{x_1,x_2,x_3\})\To \GL_7(\CC)$$be the 
monodromy representation of $\HHH.$ Then the following 
holds
\begin{itemize}
\item 
The Zariski closure of the image of $\rho$ 
coincides with $G_2(\CC).$
\item  The representation  $\rho$ (resp. the local system $\HHH$) is 
$G_2$-rigid. 
\item The Jordan canonical forms of the 
local monodromies at $x_1,x_2,x_3,\infty$ are as follows (respectively):
$$ \J(2,2,1,1,1),\,\,\J(2,2,1,1,1),\, \,
\diag(1,\zeta_3,\zeta_3,\zeta_3,\zeta_3^{-1},\zeta_3^{-1},\zeta_3^{-1}),\,\, 
\J(3,3,1)\,.$$
\end{itemize}
\end{prop}

\proof The claim on the Jordan canonical forms is obvious from the
matrices given in the Appendix. 
Since the tuples $\MC_{\zeta_3^{-1}}(\g_i),\,i=1,2,$ are contained 
in the group $\Sp_2(\CC)=\SL_2(\CC),$ their Kronecker product  $\MC_{\zeta_3^{-1}}(\g_1)\otimes \MC_{\zeta_3^{-1}}(\g_2)$
is contained in 
the orthogonal group $\SO_4(\CC)^4\leq \GL_4(\CC)^4.$
 Moreover, 
the elements of  $\MC_{\zeta_3^{-1}}(\g_1)\otimes \MC_{\zeta_3^{-1}}(\g_2)$
can easily be seen to generate an irreducible subgroup of 
$\GL_4(\CC).$ Thus, by \cite{DR00}, Cor. 3.6, or by \cite{Katz96}, Thm. 2.9.8,
the local system $\HHH$ is irreducible (i.e., the elements of $\h$ generate 
 an irreducible subgroup of $\GL_7(\CC)$). Since the elements 
of  $\MC_{\zeta_3^{-1}}(\g_1)\otimes \MC_{\zeta_3^{-1}}(\g_2)$ are
 contained in $\SO_4(\CC),$ it follows from \cite{DR00}, Cor. 5.15, that
the elements of $\h$ are contained  in the group $\SO_7(\CC).$ 
By an elementary computation (using the computer algebra system 
MAGMA \cite{magma}), one can check that the matrices 
stabilize a one-dimensional subspace of the third exterior power
$\Lambda^3(\CC^7)  .$ Thus,  by the results of \cite{Asch},
the image of $\rho$ is contained in $G_2(\CC).$  By the classification
of bireflection groups given in \cite{GuralnickSaxl}, Thm. 7.1 and Thm. 8.3, 
the Zariski closure of the 
image of $\rho$ can be seen to coincide (up to $\GL_7(\CC)$-conjugation)
with $G_2(\CC).$ By Prop. \ref{propSV},
the structure of the Jordan canonical forms of 
$h_1,\ldots,h_4$ then implies that the representation 
is $G_2$-rigid. 
\Endproof

\begin{rem}\label{remarg} \begin{enumerate} 
\item By taking 
$6$-th roots of unity instead of third roots of unity  for the definition 
of $\g_1$ and $\g_2,$ and by taking 
$${\bf s}:=\MC_{\zeta_6^{-1}}(\g_1)\otimes \MC_{\zeta_6^{-1}}(\g_2)\,,$$
one obtains another orthogonally rigid quadruple in $\SO_4(\CC),$
whose Jordan canonical forms coincide with 
those of $\MC_{\zeta_3^{-1}}(\g_1)\otimes \MC_{\zeta_3^{-1}}(\g_2).$
After application of $\g_3^{-1} \otimes \MC_{\zeta_3^{-1}}\circ 
\g_3\otimes \MC_{\zeta_3}$ one obtains an orthogonally rigid tuple $\h'$ 
in the group $\SO_7(\CC)$
which has the same tuple of Jordan canonical forms as the above 
tuple $\h\in G_2(\CC)^4.$ 
(Thus the containment of the elements of $\h$ in the group 
$G_2$ cannot be derived from the information on the Jordan canonical forms 
alone.)
\item Proof of Theorem 1 (iii): Using the tensor structure 
of the group $\SO_4(\CC),$ it 
 can be checked that any irreducible and orthogonally rigid
tuple in the group $\SO_4(\CC)$ is $\GL_4(\CC)$-conjugate either to 
the above tuple ${\bf s}$ or to 
$\MC_{\zeta_3^{-1}}(\g_1)\otimes \MC_{\zeta_3^{-1}}(\g_2).$ It follows 
then from the invertibility  of $\MC_\lambda$ (see \cite{DR00}, Thm. 5.3,
or \cite{Katz96}, 5.1.5) that any tuple in the group $\SO_7(\CC)$ which has
the same tuple of Jordan canonical forms is either $\GL_7(\CC)$-conjugate
to $\h$ or to $\h'.$  \end{enumerate}
\end{rem}

\begin{prop}\label{Prop44} {\rm (Proof of Theorem 1 (iii) and (iv))}
\begin{enumerate}
\item The local  system $\HHH$ is $G_2$-rigid 
and $\SO_7$-rigid but not 
$\GL_7$-rigid.
\item The rigidity conjecture holds for the $G_2$-rigid local system
$\HHH.$
\end{enumerate}
\end{prop}

\proof  The first claim follows from Theorem \cite{SV} considering the 
Jordan canonical forms of the monodromy tuple of $\HHH,$ given in 
 Prop. \ref{Properg1}.
The second claim follows by an iterative application 
of the K\"unneth-formula and Prop. \ref{propkat}. 
\Endproof

\begin{rem}\begin{enumerate}
\item Let $3<q$ be a natural number and $\zeta_q$ be a primitive 
$q$-th root of unity. Let $\g_1$ be as above and 
let $\g_2'=(1,\zeta_q\zeta_3,\zeta_q^{-1}\zeta_3
,\zeta_3).$ Let 
further $\g_{q,1}=(1,\zeta_q^{-1},1,\zeta_q),$
 $\g_{q,2}=(1,\zeta_q^{-1},\zeta_q,1),$
 $\g_{q,3}=(1,1,\zeta_q^{-1} \zeta_3,\zeta_3^{-1}\zeta_q)$
 and
  $\g_{q,4}=(1,1,\zeta_3^{-1},\zeta_3).$
 Using the sequence 
\begin{eqnarray}\g_{q,4} \otimes \MC_{\zeta_q \zeta_3}(\g_{q,3} \otimes( \MC_{\zeta_3^{-1}}(\g_{q,2}\otimes( \MC_{\zeta_q^{-1}}(\g_{q,1}^{-1}\otimes\quad \quad  && \nonumber 
\\ 
\quad \quad (\MC_{\zeta_3^{-1}}(\g_1)\otimes (\g_{q,1} \otimes \MC_{\zeta_q \zeta_3^{-1}}(\g_2'))))))\nonumber&&\end{eqnarray}
one obtains tuples $\h_q$ whose tuples of Jordan canonical forms 
are 
$$ (\J(2,2,1,1,1),\,\,\J(2,2,1,1,1),\, \,
\diag(1,\zeta_3,\zeta_3,\zeta_3,\zeta_3^{-1},\zeta_3^{-1},\zeta_3^{-1})\,,$$
$$ 
 \diag(\zeta_q,\zeta_q,1,1,1,\zeta_q^{-1},\zeta_q^{-1})).$$
The authors expect these tuples to have similar properties as 
the above considered tuple $\h$ but lack of computational power prevents
a proof of these statements in the general case. 
\item Assume for simplicity 
that $\{x_1,x_2,x_3,\infty\}=\{0,\pm 1 ,\infty\}$
and let 
 $$f:\PP^1\setminus \{0,\pm 1,\infty\}
\to \PP^1 \setminus \{0,1,\infty\},\, x\mapsto x^2\,.$$ 
One can show that the local system 
$\HHH$ is the pullback $\HHH=f^*{\cal F},$
where ${\cal F}$ is the $G_2$-rigid local system associated to a triple
$(f_1,f_2,f_3)\in G_2(\CC)^3$ whose Jordan canonical forms 
are $\diag(1,-\zeta_3,-\zeta_3,\zeta_3,-\zeta_3^{-1},-\zeta_3^{-1},\zeta_3^{-1}),\, \J(2,2,1,1,1), $ and an element having one Jordan block of length one 
and eigenvalue $-1,$ one Jordan block of length three
and eigenvalue $-1,$ and one Jordan block of length three
and eigenvalue $1$ (here, $f_1$ gives the monodromy at $0,$ $f_2$ gives the monodromy at $1,$ and $f_3$ gives the monodromy at $\infty$). 
\end{enumerate}
\end{rem}

\section{Appendix: The explicit monodromy tuple of $\HHH$}\label{App}
The monodromy tuple of the local system $\HHH\in \LS(\PP^1\setminus
\{x_1,x_2,x_3\})$ is given by 
$(h_1,\ldots,h_4),$ where $h_4=(h_1h_2h_3)^{-1}$ and 
$h_1,h_2,h_3\in \GL_7(\CC)$ are as follows:

\begin{small}
$$ h_1=\left(\begin{array}{ccccccc}
     1& -3& \zeta_3 - 1& 0& \zeta_3 - 4 &0& 2\zeta_3 + 4\\
     0& 3\zeta_3 + 1& 2\zeta_3 + 1& 0& 2\zeta_3 + 1& -2\zeta_3 - 1& 0\\
     0& -3\zeta_3 &-2\zeta_3& 0 &-2\zeta_3 - 1& 2\zeta_3 + 1 &0\\
     0& 3\zeta_3 + 3& \zeta_3 + 2 &1& \zeta_3 + 2& -\zeta_3 - 2& 0\\
     0 &3\zeta_3 + 6& 3& 0& 4& -3& 0\\
     0& 3\zeta_3 + 3& \zeta_3 + 2& 0& \zeta_3 + 2& -\zeta_3 - 1& 0\\
     0& 6& -2\zeta_3 + 2 &0 &-2\zeta_3 + 2& 2\zeta_3 - 2& 1
\end{array}\right)$$
$$ h_2=\left(\begin{array}{ccccccc}
   1& 0& 0& 0& 0& 0& 0\\
   \zeta_3 - 1& 1& 0& 2\zeta_3 + 1& 0& 0& 0\\
     3& 0 &1 &-2\zeta_3 - 1& -3& 0& 2\zeta_3 + 4\\
     0& 0 &0 &1& 0& 0& 0\\
     0 &0& 0& 0& 1& 0& 0\\
     0& 0& 0 &0& 0& 1 &0\\
     0 & 0& 0 &0& 0& 0 &1\end{array}\right),$$ $$
 h_3=\left(\begin{array}{ccccccc}
     \zeta_3 &0 &0 &0& 0& 0& 0\\
     0& \zeta_3 &0& 0& 0& 0& 0\\
     0& 0& \zeta_3& 0& 0& 0& 0\\
     \zeta_3 + 2& 0& 0& -\zeta_3 - 1& 0& 0& 0\\
     0& \zeta_3 + 2& 0& 0& -\zeta_3 - 1& 0& 0\\
     0 &3\zeta_3 + 3& \zeta_3 + 2& 0& 0& -\zeta_3 - 1& 0\\
     0& 0& 0& 0& \zeta_3 - 1& 0& 1\end{array}\right)\,.$$

\end{small}

\bibliographystyle{plain} \bibliography{p}

Michael Dettweiler

IWR, Universit\"at Heidelberg,

INF 368

69121 Heidelberg, Deutschland

e-mail: michael.dettweiler@iwr.uni-heidelberg.de\\

Stefan Reiter

Technische Universit\"at Darmstadt

Fachbereich Mathematik AG 2

Schlo\ss gartenstr. 7

64289 Darmstadt, Deutschland

email: reiter@mathematik.tu-darmstadt.de

\end{document}